\documentclass[11pt,twoside]{article}
\usepackage[T1]{fontenc}
\usepackage{multicol}
\usepackage[english]{babel}
\usepackage{fancybox}
\usepackage{graphicx}
\usepackage[utf8]{inputenc}
\usepackage{enumerate}
\usepackage[colorlinks=true]{hyperref}
\usepackage{amsfonts,amsmath,subfigure,amssymb,amsthm,graphicx,epsfig,float,url, bm}
\usepackage{color}
\usepackage{tikz}
\usepackage{pgf,tikz}
\usepackage{fullpage}
\usepackage{mathrsfs}
\usepackage{dsfont}
\usepackage{ bbold }
\usepackage{comment}
\usepackage{stmaryrd}
\usepackage[toc,page]{appendix}

\numberwithin{equation}{section}
\usepackage{titlesec}
\titleformat*{\section}{\Large\bfseries\raggedright}
\titleformat*{\subsection}{\large\bfseries\raggedright}

\theoremstyle{plain}
\newtheorem{thrm}{Theorem}[section]
\newtheorem{lmm}[thrm]{Lemma}
\newtheorem{crllr}[thrm]{Corollary}

\newtheorem{hyp}{{Assumption}}[section]
\theoremstyle{definition}

\newtheorem{rmrk}[thrm]{Remark}

\theoremstyle{plain}

  \def\xR{\mathbb{R}}
  \def\xN{\mathbb{N}}

\providecommand{\keywords}[1]{\textbf{{Keywords:}} #1}
\providecommand{\subjclass}[1]{\textbf{{MSC:}} #1}
 \usepackage[colorlinks=true]{hyperref}

 \numberwithin{equation}{section}
 
 \title{\textbf{A Silverman--Meadows--type criterion for the null controllability of time-dependent coupled heat equations with non-diagonalizable diffusion matrices}}

\author{
Kala Agbo bidi\thanks{XLIM, Universit\'e de Limoges, CNRS, Limoges, France.}
\thanks{IMB, Universit\'e de Bourgogne Europe, CNRS, Dijon, France.}
\and
Thibault Liard\footnotemark[1]
\and
Pierre Lissy\thanks{CERMICS, CNRS, ENPC, Institut Polytechnique de Paris, Marne-la-Vall\'ee, France.}
}

\date{\empty}

\begin{document}
\maketitle

\begin{abstract} This paper  is devoted to a uniform spectral Silverman--Meadows rank condition for the null--controllability of coupled systems of heat equations with time-dependent diffusion, coupling and control matrices. The proof relies on an adaptation of the Lebeau--Robbiano strategy combined with uniform finite-dimensional observability estimates associated with each Laplace eigenmode. \end{abstract}



\maketitle

\section{Introduction}

The controllability and observability of finite--dimensional linear systems are
classically characterized by rank conditions and Gramian criteria. For autonomous
systems, the fundamental algebraic reference is the Kalman rank condition
\cite{KalmanHoNarendra1963}. For time-dependent systems, controllability on a time interval can be characterized by the non-degeneracy of a suitable
controllability Gramian (see \cite{KalmanHoNarendra1963}). In parallel, more local
algebraic criteria for time-dependent systems, involving successive derivatives of the control matrix along the dynamics, were introduced by Silverman and Meadows \cite{SilvermanMeadows1967}. 
In the present paper, our goal is to partially extend this algebraic viewpoint
to coupled systems of partial differential equations controlled or observed
through a reduced number of components. 
Extending these algebraic criteria to infinite-dimensional problems offers several significant advantages. Once the scalar observability or spectral estimates are available, the effect of the coupling and control matrices
can be analyzed through explicit algebraic conditions, in contrast to approaches based on delicate analytical estimates such as multiplier methods, Carleman inequalities, or microlocal  and semiclassical analysis. 
A mode-by-mode Gramian criterion characterizes the controllability of each eigenmode, but gives no information on its cost. The method used in this paper, a Lebeau–Robbiano-type strategy, requires mode-by-mode observability estimates on arbitrarily small time intervals, with constants explicit in the length of the interval and polynomial in the frequency; spectral Silverman–Meadows conditions are more adapted to provide these quantitative bounds.

The literature on the controllability and observability of coupled PDE systems
is very large. We do not attempt to give here
an exhaustive account of this literature. Rather, we focus on the part of it
where controllability or observability is governed by algebraic conditions on
the coupling and control matrices, in the spirit of Kalman and
Silverman--Meadows criteria.

For parabolic systems, algebraic controllability criteria have been investigated
in several directions. Silverman--Meadows and spectral Kalman type conditions
were obtained in
\cite{AmmarKhodjaBenabdallahDupaixGonzalez2009DEA,
AmmarKhodjaBenabdallahDupaixGonzalez2009JEE}, extended in \cite{Olive2012} to
controls supported in distinct subdomains, and further developed for boundary
control and non-diagonalizable systems in
\cite{AmmarKhodjaBenabdallahGonzalezDeTeresa2011,
FernandezCaraGonzalezBurgosDeTeresa2015,DuprezGonzalezBurgosSouza2023}; see also the survey 
\cite{AmmarKhodjaBenabdallahGonzalezDeTeresa2011Survey}. In
\cite{LissyZuazua2019}, a general algebraic approach was developed to transfer
scalar observability inequalities to coupled PDE systems. In the specific case
of heat equations with constant coefficients, this also led to a sharp spectral
Kalman characterization in the natural $L^2$ setting, allowing
non-diagonalizable diffusion matrices and observations on possibly distinct
subdomains. A different (but related) algebraic viewpoint, based on the algebraic solvability of underdetermined systems rather than on rank conditions, goes back to \cite{CoronLissy2014} and was developed in
\cite{DuprezLissy2016,DuprezLissy2018,
SteevesGharesifardMansouri2019a,SteevesGharesifardMansouri2019b}. We also
mention recent Kalman-type results for coupled Stokes systems
\cite{TakahashiDeTeresaWuZhang2024,LeBalchDeTeresa2025}.
Mixed systems, involving components with different qualitative behaviours, form another important class where algebraic and analytic effects interact strongly. The parabolic--transport systems studied in
\cite{BeauchardKoenigLeBalch2020} provide a first systematic framework in this direction. This was further developed in \cite{KoenigLissy2026}, where a sharp spectral Kalman rank condition is obtained for underactuated
constant-coefficient parabolic--transport systems in large time, together with negative results in small time and for rough initial data.

The present paper continues this line of research. Its purpose is to derive
algebraic criteria for the controllability or observability of coupled PDE
systems through a reduced number of controls or observations. We consider  class of time-dependent coupled  linear parabolic equations where each control may act on a different subdomain of the spatial domain. Our main contribution is the use of the direct  Lebeau--Robbiano strategy introduced by Miller \cite{MR2679651} to prove the null-controllability under the uniform Silverman--Meadows rank condition. More precisely, after expanding the solution in an orthonormal
basis of Dirichlet Laplacian eigenfunctions, each spectral mode  is associated with a finite-dimensional time-dependent control system. We introduce a family of Silverman--Meadows type matrices depending both on time and on Laplacian eigenvalues and prove that if these matrices satisfy a rank condition uniformly with respect to both the time variable and the spectral parameter, then the original parabolic system is null-controllable. The main difficulty in proving the result  lies in the fact that, unlike the framework considered in \cite{MR2679651}, we cannot exploit the translation property of the semigroup, since the underlying operator is non-autonomous. Instead, the analysis must be carried out through the associated evolution resolvents, which makes the argument substantially more involved. In this context, a uniform Silverman--Meadows condition plays a crucial role. Overall, our contribution extends Miller's direct Lebeau--Robbiano method to linear parabolic systems with time-dependent diffusion, coupling, and control matrices.

The paper is organized as follows: we first establish null controllability of the parabolic system with time-dependent coefficients under the uniform Silverman–Meadows condition (Assumption~\ref{hyp:USM}), and then derive, under additional assumptions, a pointwise condition implying this assumption.

\section{Algebraic criteria for the controllability of coupled linear  parabolic equations}
\label{sec:heat}
\subsection{Presentation of the problem and main result}
Let $\Omega$ be a smooth bounded domain of
$\mathbb R^N$, $N\in\mathbb N^*$, and let $T>0$. We consider
\begin{equation}\label{eq:intro_heat}
\left\{
\begin{array}{lll}
\partial_t Y &= D(t)\Delta Y + A(t)Y + B(t)W(x)U,
    & (t,x)\in (0,T)\times\Omega,\\
Y&=0, & (t,x)\in (0,T)\times\partial\Omega,\\
Y(0,\cdot)&=Y^0, & x\in\Omega,
\end{array}
\right.
\end{equation}
where $Y^0\in L^2(\Omega)^n$, $U\in L^2((0,T)\times\Omega)^m$, and $1\leq m\leq n$ two integers.
The matrix $W$ describes the localization of the different controls and may
involve distinct open subsets of $\Omega$. More precisely, we consider
\begin{equation}
    W(x)=
\begin{pmatrix}
\mathbf 1_{\omega_1}(x) & 0 & \cdots & 0\\
0 & \mathbf 1_{\omega_2}(x) & \ddots & \vdots\\
\vdots & \ddots & \ddots & 0\\
0 & \cdots & 0 & \mathbf 1_{\omega_m}(x)
\end{pmatrix}\in M_m(\xR),
\end{equation}
where $\omega_1,\ldots,\omega_m$ are nonempty open subsets of $\Omega$ (that might be distinct). We will assume in what follows
that 
\begin{equation}\label{smoothADB}
    A\in C^{\infty}(\xR_+;M_{n}(\xR)), D\in C^{\infty}(\xR_+;M_{n}(\xR)), B\in C^{\infty}(\xR_+;M_{n,m}(\xR)),
\end{equation}
and that $D(t)$ is uniformly
coercive, namely, there exists $\alpha>0$ such that
\begin{equation}\label{coercivityofD}
\langle D(t)X,X\rangle_{\mathbb R^n}\geq \alpha \|X\|^2,
\;\;
\forall X\in\mathbb R^n,\; \forall t\in[0,T].
\end{equation}

Let $V:=H^1_0(\Omega)^n$, $H:=L^2(\Omega)^n$ and
\[
\mathcal{W}(0,T):=\bigl\{v\in L^2(0,T;V):\ \partial_tv\in L^2(0,T;V')\bigr\},
\]
which embeds continuously into $C([0,T];H)$
\cite[Chap.~XVIII, \S1, Theorem~1]{DL92}.
For $t\in[0,T]$ and $u,v\in V$, set
\[
a(t;u,v):=\sum_{i=1}^N\int_\Omega\bigl\langle D(t)\partial_iu,\partial_iv\bigr\rangle\,dx
-\int_\Omega\bigl\langle A(t)u,v\bigr\rangle\,dx
+\lambda_A\int_\Omega\langle u,v\rangle\,dx,
\qquad
\lambda_A:=\|A\|_{L^\infty(0,T)} .
\]
By \eqref{smoothADB} and \eqref{coercivityofD} applied to $X=\partial_iu(x)$,
$a$ is bounded on $V\times V$ and
\[
a(t;u,u)\ \ge\ \alpha\,\|\nabla u\|_{L^2(\Omega)}^2
\qquad\forall u\in V,\ \forall t\in[0,T].
\]
 Hence, by Lions' theorem
\cite[Chap.~XVIII, \S3, Theorem~2]{DL92}, for every $Y^0\in H$ and
$U\in L^2((0,T)\times\Omega)^m$ there is a unique $\tilde Y\in \mathcal{W}(0,T)$ such that
$\tilde Y(0)=Y^0$ and
\[
\bigl\langle\partial_t\tilde Y(t),v\bigr\rangle_{V',V}+a(t;\tilde Y(t),v)
=e^{-\lambda_A t}\int_\Omega\bigl\langle B(t)W(x)U(t),v\bigr\rangle\,dx
\qquad\forall v\in V,\ \text{for \textit{a.e.} } t\in(0,T).
\]
Then $Y:=e^{\lambda_A t}\tilde Y\in \mathcal{W}(0,T)$ is the unique weak solution of
\eqref{eq:intro_heat}, and
\[
\|Y\|_{C([0,T];H)}+\|Y\|_{\mathcal{W}(0,T)}
\ \le\ C\bigl(\|Y^0\|_H+\|U\|_{L^2((0,T)\times\Omega)^m}\bigr).
\]
Since $\langle D^*X,X\rangle=\langle X,DX\rangle$, the adjoint system (that will be introduced in Section~\ref{proofth1}) can be handled in the same way after the change of variable
$t\mapsto T-t$.

Before stating the main theorem, we introduce some notation that will be used throughout the paper.  Let $(\lambda_k)_{k\geq1}$ be the non-decreasing sequence of eigenvalues of
$-\Delta$ with Dirichlet boundary condition, repeated according to their
multiplicity, and let $(e_k)_{k\geq1}$ be an associated orthonormal basis of
$L^2(\Omega)$ made of eigenfunctions, $e_k\in H_0^1(\Omega)$. For $\lambda\in\xR$,
set
\begin{equation}\label{def:A_lambda}
    A_{\lambda}(t):=A(t)-\lambda D(t),
\end{equation}
and let us define, by induction on $i$, a sequence of maps $B_i^{\lambda}$ in
the following way:
\begin{equation}\label{def:kalman}
    B_0^{\lambda}(t):=B(t),\;\;
B_i^{\lambda}(t):=\dot B_{i-1}^{\lambda}(t)-A_{\lambda}(t)B_{i-1}^{\lambda}(t),
\; i\geq 1.
\end{equation}
Note that, by \eqref{smoothADB},
\begin{equation}\label{eq:smoothB_i}
    B_i^{\lambda} \in C^\infty\bigl([0,T];{M}_{n,m}(\mathbb{R})\bigr),
    \qquad\forall\lambda\in\xR,\ \forall i\in\xN .
\end{equation}
For $p\in\xN$ and $\lambda\in\xR$, we set the following matrix function
\begin{equation}\label{def:kalmanMatrix}
    \mathcal K_{\lambda,p}(t):=\bigl[
B_0^\lambda(t),B_1^\lambda(t),\ldots,B_p^\lambda(t)
\bigr]\in  M_{n,(p+1)m}(\xR).
\end{equation}
and the scalar function 
\begin{equation}\label{def:Pp}
P_p(\lambda,t)
:=
\det\bigl(
\mathcal K_{\lambda,p}(t)\mathcal K_{\lambda,p}(t)^T
\bigr).
\end{equation}
 \ \\
Note that we have $\operatorname{rank}(\mathcal K_{\lambda,p}(t)\mathcal K_{\lambda,p}(t)^T) =
\operatorname{rank}\mathcal K_{\lambda,p}(t)$. Thus,
\begin{equation}\label{eq:equi_rank}
    P_p(\lambda,t)>0
\;\Longleftrightarrow\;
\operatorname{rank}\mathcal K_{\lambda,p}(t)=n.
\end{equation}
\begin{rmrk}\label{rem:polynomial}
For every $i\in\{0,\ldots,p\}$, from \eqref{def:A_lambda} and \eqref{def:kalman}, there exist matrices $V_j^i\in C^\infty\bigl([0,T];{M}_{n,m}(\mathbb{R})\bigr)$ such that 
\begin{equation*}
    B_i^\lambda(t)=\lambda^i D^i(t) B(t)+\sum_{j=0}^{i-1}V_j^i(t)\lambda^j,\;\forall\;\lambda\in\xR,
\end{equation*}
In particular,
\begin{equation*}
    \deg_\lambda B_i^\lambda(t)\leq i.
\end{equation*}
Since each entry of $\mathcal K_{\lambda,p}(t)\mathcal K_{\lambda,p}(t)^T$ is a polynomial in $\lambda$ of degree at most $2p$, we deduce that $ \lambda \mapsto P_p(\lambda,t)$ is polynomial and 
\begin{equation}
    \deg_\lambda P_p(\lambda,t)\leq 2pn.
\end{equation}
\end{rmrk}\label{rem:degree}
Our main result relies on the following uniform spectral Silverman–Meadows condition.
\begin{hyp}[Uniform spectral Silverman--Meadows condition]
\label{hyp:USM}
There exist $p\in\xN$ and a nontrivial closed interval $[a,b]\subset[0,T]$ such that
\begin{equation}\label{hyp:USM:rank}
\operatorname{rank}\mathcal K_{\lambda_k,p}(t)=n,
\;
\forall t\in[a,b],\; \forall k\in\xN^*,
\end{equation}
or, equivalently, $  \operatorname{Span}
\left\{
B_i^{\lambda_k}(t)X;\ X\in\xR^m,\ i=0,\ldots,p
\right\}
=\xR^n,
\;
\forall t\in[a,b],\; \forall k\in\xN^*$. 
\end{hyp}
The main result is the following.
\begin{thrm}\label{main:thm}
Assume that Assumption~\ref{hyp:USM} holds. Then
System~\eqref{eq:intro_heat} is null-controllable at time $T$: for any $Y^0\in L^2(\Omega)^n$, there exists $U\in L^2((0,T)\times \Omega)^m$ such that $Y(T,\cdot)=0$ in $L^2(\Omega)^n$.
\end{thrm}

Let us explain precisely what our result adds to the existing literature.  Compared with the constant-coefficient parabolic result of Lissy and
Zuazua~\cite{LissyZuazua2019}, our theorem allows the diffusion, coupling and control
matrices to depend on time while retaining non-diagonalizable diffusion
matrices. Time-dependent coupling and control matrices
were already treated by Ammar Khodja et al.~\cite{AmmarKhodjaBenabdallahDupaixGonzalez2009DEA} with a common scalar
parabolic operator, and by Olive~\cite{Olive2012} with time-independent coefficients, $D=I_n$ and controls supported
in distinct open subsets. 

\begin{rmrk}
Theorem~\ref{main:thm} is proved under the uniform condition
Assumption~\ref{hyp:USM}. Unlike the finite-dimensional setting of
\cite[Theorem~1.18]{MR2302744}, where the Silverman--Meadows condition at a
single time is sufficient, the pointwise spectral Silverman--Meadows condition,
namely the rank condition at a single time $\bar t$ for the whole family of
matrices $\mathcal K_{\lambda_k,p}$, $k\geq 1$ (Assumption~\ref{hyp:pointwise}
below), does not imply in general a rank condition that is uniform both in
time, on some nontrivial interval, and in the frequency $k$; see
Remark~\ref{rem:notpointwise} for an example with smooth coefficients.
Sufficient conditions under which the pointwise condition does imply
Assumption~\ref{hyp:USM} are given in Section~\ref{sec:sufficient}. Outside
these cases, our method does not apply.  Whether the pointwise
condition alone implies null-controllability for smooth coefficients, or
whether there exist pathological examples for which null-controllability
fails, is an open question that we were not able to address.
\end{rmrk}
\begin{rmrk}\label{rem:strategy}
The proof of Theorem~\ref{main:thm} follows the direct Lebeau--Robbiano
strategy of Miller~\cite{MR2679651}. Recall that the Lebeau--Robbiano method
\cite{95LR} (see also \cite[pp.~19--21]{olive2013} for detailed explanations) proves the
null-controllability of a parabolic equation $\dot Z=\mathcal AZ$ on a Hilbert
space $H$, where $-\mathcal A$ is a nonnegative self-adjoint operator
with compact resolvent, from two ingredients: a spectral inequality
\begin{equation}\label{ineq:stationnaire}
\|u\|_H^2\le c_0e^{c\lambda^\alpha}\|Cu\|_Y^2,\qquad u\in H_\lambda,
\end{equation}
where $H_\lambda$ is spanned by the eigenfunctions with eigenvalue at most
$\lambda$, $C$ is the observation operator and $\alpha\in(0,1)$, and the
dissipation $$\|e^{t\mathcal A}u\|_H\le e^{-t\lambda}\|u\|_H, \qquad u\perp H_\lambda$$
In its original form, the method alternates, on smaller
and smaller time intervals, a control phase on $H_\lambda$ and a free
dissipation phase, with an increasing threshold $\lambda$; the choice of the
threshold at each step balances the spectral factor $e^{c\lambda^\alpha}$
against the damping factor. Miller's direct strategy \cite[Lemma~2.1]{MR2679651}
avoids the construction of controls altogether: it combines
\eqref{ineq:stationnaire}, at a single well-chosen frequency scale depending
on the length $\varepsilon$ of the time interval, with a purely temporal
telescoping argument. This
yields a short proof of the observability inequality, and it applies to broader abstract settings
(non-self-adjoint generators, non-discrete spectrum), as highlighted in \cite[Theorem 2.1]{beauchard2018null}.

Here, we apply this mechanism to the non-autonomous system \eqref{eq:obs}: the
low-frequency estimate is provided mode by mode by Lemma~\ref{observabilityineq},
whose constant is explicit in $\varepsilon$ and polynomial in $\lambda_k$, and
the spectral inequality \eqref{ineq:stationnaire} is the Lebeau--Robbiano
inequality for the Dirichlet Laplacian ($\alpha=\tfrac12$). Two adaptations
are needed. First, the absence of a semigroup structure rules out the time
translation used by Miller to propagate the local estimate; this is why
Lemma~\ref{MillerLemma} is stated for every subinterval
$[c,c+\varepsilon]\subset[\tilde a,\tilde b]$, with $f$ and $q$ independent of $c$, and why
Lemma~\ref{observabilityineq} is applied on the terminal window
$[c+\varepsilon/2,c+\varepsilon]$ only, the initial half $[c,c+\varepsilon/2]$
letting the dissipation of high frequencies act. Second, the telescoping is
carried out towards the endpoint $c$ of the subinterval rather than towards
$0$, which gives \eqref{ine:obscost}.
\end{rmrk}
\subsection{Proof of Theorem~\ref{main:thm} }
\label{proofth1}
Before giving the proof of Theorem~\ref{main:thm}, we are going to prove some intermediate results.

Our first lemma actually proves a quantitative version of
\eqref{hyp:USM:rank}, namely a uniform lower bound on the Gram determinant
$P_p(\lambda_k,t)$, that is automatically
available on a subinterval. Recall that the recursion \eqref{def:kalman} defines
$B_i^{\lambda}(t)$ for every real $\lambda$, so that, by
Remark~\ref{rem:degree}, $P_p(\lambda,t)=\sum_{s=0}^{2pn}\pi_s(t)\lambda^s$ with
$\pi_s\in C^\infty([0,T])$, and $P_p(\lambda,t)\ge0$ for all
$(\lambda,t)\in\xR\times[0,T]$, as the determinant of a Gram matrix.

\begin{lmm}\label{lem:uniformGram}
Assume that Assumption~\ref{hyp:USM} holds. Then there exist a
constant $c_G>0$ and a nontrivial closed
interval $[a',b']\subset[a,b]$ such that
\begin{equation}\label{hyp:USM:gram}
P_p(\lambda_k,t)\ \ge\ c_G,
\qquad
\forall t\in[a',b'],\ \forall k\in\xN^*.
\end{equation}
\end{lmm}

\begin{proof}
Since $P_p(\lambda_1,\cdot)>0$ on $[a,b]$ by \eqref{hyp:USM:rank} and
\eqref{eq:equi_rank}, the coefficients $\pi_s$ do not all vanish identically on
$[a,b]$, and we may set
\[
d:=\max\bigl\{s\in\{0,\ldots,2pn\}:\ \pi_s\not\equiv0 \text{ on } [a,b]\bigr\},
\qquad
O:=\{t\in[a,b]:\ \pi_d(t)\neq0\}.
\]
Since $\pi_s\in C^\infty([0,T])$, $O$ is a nonempty relatively open subset of $[a,b]$, and
$P_p(\lambda,t)=\sum_{s\le d} \pi_s(t)\lambda^s$ for $t\in[a,b]$. Since
$P_p(\lambda,t)\ge0$ for all $\lambda\in\xR$, dividing by $\lambda^d$ and letting $\lambda\to+\infty$ gives
$\pi_d(t)\ge0$ for every $t\in[a,b]$; hence $\pi_d>0$ on $O$. Choose a nontrivial
closed interval $[a',b']\subset O$ and set
\[
c_0:=\min_{[a',b']}\pi_d>0,
\qquad
M_0:=\max_{0\le s<d}\ \max_{[a',b']}|\pi_s|,
\]
with $M_0:=0$ if $d=0$. For $t\in[a',b']$ and $\lambda\ge1$, using
$\sum_{s<d}\lambda^s\le d\lambda^{d-1}$,
\begin{equation*}
    P_p(\lambda,t)\ \ge\ c_0\lambda^d-dM_0\lambda^{d-1}
=\lambda^d\Bigl(c_0-\frac{dM_0}{\lambda}\Bigr)\ \ge\ \frac{c_0}{2}\lambda^d\ \ge\ \frac{c_0}{2},
\end{equation*}
$ \text{whenever }\ \lambda\ge\lambda_0:=\max\Bigl(1,\frac{2dM_0}{c_0}\Bigr).$
Since $\lambda_k\to+\infty$, the set $\{k\in\xN^*:\ \lambda_k<\lambda_0\}$ is
finite. For each such $k$, the function $P_p(\lambda_k,\cdot)$ is continuous
and, by \eqref{eq:equi_rank} and \eqref{hyp:USM:rank}, positive on the compact
interval $[a',b']\subset[a,b]$, hence bounded from below by some $c_k>0$.
Setting $c_G:=\min\bigl(c_0/2,\ \min_{\lambda_k<\lambda_0}c_k\bigr)>0$ yields
\eqref{hyp:USM:gram}.
\end{proof}

\begin{rmrk}
Even under Assumption 2.1, the lower bound (2.14) cannot, in general, be obtained on $[a,b]$. Indeed, consider  $[0,T]=[a,b]=[0,1]$, $n=m=1$, $D=1$, $A=0$, $B(t)=t-\tfrac12$ and $p=1$. From \eqref{def:A_lambda}, \eqref{def:kalman}, \eqref{def:kalmanMatrix} and \eqref{def:Pp}, one
has, for any $(\lambda,t) \in \xR \times [0,1]$,
$$P_1(\lambda,t)=(t-\tfrac12)^2+(1+\lambda(t-\tfrac12))^2>0.$$
Thus, from \eqref{eq:equi_rank}, Assumption 2.1 holds. Let $\lambda \in \xR$, we have $\min_{t\in[0,1]}P_1(\lambda,t)=(1+\lambda^2)^{-1}$ and the unique minimizer is 
\[
t_\lambda=\frac12-\frac{\lambda}{1+\lambda^2} \in [0,1].
\]
Since $\lambda_k\to+\infty$ as $k\to+\infty$, we obtain
$$
P_1(\lambda_k,t_{\lambda_k})
=\frac{1}{1+\lambda_k^2}\longrightarrow0 \text{ and }
t_{\lambda_k}\longrightarrow\frac12.$$
Thus, \eqref{hyp:USM:gram} fails on $[a,b]=[0,1]$.
\end{rmrk}

In view of Lemma~\ref{lem:uniformGram}, shrinking $[a,b]$ if necessary, we will always assume in the following 
that \eqref{hyp:USM:gram} holds on $[a,b]$.

Now, we will develop some duality theory needed in our proof.
Let  $\varphi_T\in L^2(\Omega)^n$. The adjoint system associated with
\eqref{eq:intro_heat} is the backward system
\begin{equation*}
\left\{
\begin{array}{lll}
-\partial_t\varphi
&=
D^*(t)\Delta\varphi+A^*(t)\varphi,
& (t,x)\in(0,T)\times\Omega,\\
\varphi&=0,
& (t,x)\in(0,T)\times\partial\Omega,\\
\varphi(T)&=\varphi_T,
& x\in\Omega,
\end{array} \right .
\end{equation*}
which enjoys the same well-posedness properties as the original system. We perform the time reversal $Z(t)=\varphi(T-t)$  to follow more directly Miller's argument of \cite{MR2679651} {and set}
{
\begin{equation}\label{smooth_tilde_ABD}
    \widetilde A(t)=A(T-t),
\;
\widetilde D(t)=D(T-t),
\;
\widetilde B(t)=B(T-t).
\end{equation}
}
Then $Z$ solves
{
\begin{equation}\label{eq:obs}
\left\{
\begin{array}{lll}
\partial_t Z
&=
\widetilde D^*(t)\Delta Z+\widetilde A^*(t)Z,
& (t,x)\in(0,T)\times\Omega,\\
Z&=0,
& (t,x)\in(0,T)\times\partial\Omega,\\
Z(0)&=\varphi_T,
& x\in\Omega,
\end{array} \right .
\end{equation}
with
$\widetilde A^*\in C^{\infty}([0,T];M_n(\xR))$ and
$\widetilde D^*\in C^{\infty}([0,T];M_n(\xR))$.
The time-reversed observation is $\bar B^*(t)Z(t)$, where
\begin{gather}\label{def:observ-operator}
    (\bar B^*(t)v)(x) := W^*(x)\widetilde B^*(t)v(x),
    \qquad v\in L^2(\Omega)^n.
\end{gather}
The operator $\bar B^*(t)$ maps $L^2(\Omega)^n$ into the observation space
$L^2(\Omega)^m$.
}
For the reversed coefficients, we set
\begin{equation*}
    \widetilde A_\lambda(t)
=
\widetilde A(t)-\lambda\widetilde D(t)
=A_\lambda(T-t).
\end{equation*}
For a general parameter $\lambda$, we introduce the
auxiliary finite-dimensional control system
\begin{equation}\label{eq:dualcontrol}
    \dot y(t)=-\widetilde A_{\lambda}(t)\,y(t)+\widetilde B(t)\,u(t) ,
\end{equation}
with $y(t)\in\xR^n$ and $u(t)\in\xR^m$. Its adjoint equation in the controllability--observability duality is
$\dot z=\widetilde A_\lambda^*(t)z,$
with the finite-dimensional output $\widetilde B^*(t)z\in\xR^m$.

The Silverman--Meadows matrices associated with \eqref{eq:dualcontrol} are\begin{equation}\label{eq:SMdual}
    C_0^\lambda(t)=\widetilde B(t),
\;\;
C_i^\lambda(t)
=
\frac{d}{dt}C_{i-1}^\lambda(t)
+
\widetilde A_\lambda(t)\,C_{i-1}^\lambda(t),
\;\; i\geq 1 .
\end{equation}
A straightforward induction, using \eqref{def:kalman}, \eqref{smooth_tilde_ABD} and
$\widetilde A_\lambda(t)=A_\lambda(T-t)$, gives
\begin{equation*}
    C_i^\lambda(t)=(-1)^{i}\,B_i^\lambda(T-t),\qquad i\ge0 .
\end{equation*}
Consequently, $C_i^{\lambda_k}\in C^\infty([0,T];M_{n,m}(\xR))$ for
all $i,k$, and, for every $t\in[0,T]$,
\begin{equation*}
    \operatorname{rank}\big[C_0^\lambda(t),\ldots,C_p^\lambda(t)\big]
=
\operatorname{rank}\mathcal K_{\lambda,p}(T-t),
\qquad
\sum_{i=0}^pC_i^\lambda(t)\,C_i^\lambda(t)^*
=
\sum_{i=0}^pB_i^\lambda(T-t)\,(B_i^\lambda(T-t))^*,
\end{equation*}
so that the Gram determinant is unchanged:
$\det\bigl(\sum_{i\le p}C_i^{\lambda}C_i^{\lambda*}\bigr)(t)=P_p(\lambda,T-t)$.
In particular, if Assumption~\ref{hyp:USM} and \eqref{hyp:USM:gram} hold for
the original coefficients on $[a,b]$, then
\begin{equation}\label{eq:USMdual}
    \operatorname{rank}\big[C_0^{\lambda_k}(t),\ldots,C_p^{\lambda_k}(t)\big]=n
\;\;\text{and}\;\;
\det\Bigl(\sum_{i=0}^pC_i^{\lambda_k}(t)C_i^{\lambda_k}(t)^*\Bigr)\ge c_G,
\qquad\forall t\in[T-b,T-a],\ \forall k\in\xN^*,
\end{equation}
with the same integer $p$ and the same constant $c_G$.

{
Set
\begin{equation} \label{def:atilde}
\widetilde a:=T-b,\qquad\widetilde b:=T-a.
\end{equation}
We now return to \eqref{eq:obs}, observed through $\bar B^*(t)Z(t)$
as defined in \eqref{def:observ-operator}.
By well-posedness and continuity, $\bar B^*$ is an admissible observation operator for the
evolution family associated with $\widetilde D^*(t)\Delta+\widetilde A^*(t)$,
in the sense that for the fixed horizon $T$, there exists a constant $C_T>0$ such that
\begin{equation}\label{eq:admissibility}
    \int_0^T
    \left\|
        W^*\widetilde B^*(t)Z(t)
    \right\|_{L^2(\Omega)^m}^2\,dt
    \leq
    C_T
    \|Z^0\|_{L^2(\Omega)^n}^2,
    \qquad Z^0\in L^2(\Omega)^n.
\end{equation}
}

By a standard duality argument (see, \textit{e.g.}, \cite[Theorem 2.44, Page 56]{MR2302744}),  System \eqref{eq:intro_heat} is null-controllable at time $T$ if and only if  System \eqref{eq:obs} is observable in time $T$ in the sense that there exists $C>0$ such that, for every $Z^0\in L^2(\Omega)^n$, the solution $Z$ of \eqref{eq:obs} verifies
 \begin{equation} \label{ine:obs}
{\Vert Z(T) \Vert_{L^2(\Omega)^n}^2}
 \leq C \int_0^T \Vert W^*{\widetilde B}^*(t) Z(t) \Vert_{L^2(\Omega)^m}^2 dt.
 \end{equation}
Hence, we will focus on proving \eqref{ine:obs}.  We decompose $Z^0$ into a sum of eigenvectors $(e_k)_{k\in \xN^*}$: 
$$ Z^0= \sum_{k=1}^{\infty} Z_k^0 e_k(x),$$
with $Z_k^0\in \xR^n$. The solution $Z$ of \eqref{eq:obs} can be rewritten as $Z:(t,x)\mapsto \sum_{k=1}^{\infty} Z_k(t) e_k(x)$ with $Z_k$ solution of 
\begin{equation} \label{eq:decomposition}
 \left\{\begin{array}{l}
\dot Z_k(t)={\widetilde A}_{\lambda_k}^*(t)Z_k(t),\\
Z_k(0)=Z^0_k.
\end{array}\right.
 \end{equation}
 
 Let us introduce some useful notation. For $J\in\mathbb N^*$, define the two spaces, corresponding to low and high frequency spaces of functions with frequency cut at $J$:
\begin{equation}
\begin{aligned}\label{ejfj}
E_J&=
\left\{
\sum_{k=1}^{J} a_k e_k(x)\ ;\ a_k\in\mathbb R^n
\right\},
\\
F_J&=
\left\{
\sum_{k=J+1}^{+\infty} a_k e_k(x)\ ;\ a_k\in\mathbb R^n,\,
\sum_{k=J+1}^{+\infty}\|a_k\|^2<+\infty
\right\}.
\end{aligned}
\end{equation}
Note that $E_J$ and $F_J$ form an orthogonal decomposition of $L^2(\Omega)^n$, \textit{i.e.} 
\begin{equation}\label{def:direct_sum}
    L^2(\Omega)^n=E_J\oplus_{\perp} F_J,
\end{equation}
and that both are stable under the (non-autonomous) evolution family generated by $\mathcal{A}(t):={\widetilde D}^*(t)\Delta+{\widetilde A}^*(t)\,\mathrm{Id}$
denoted by $(U(t,s))_{{0\leq s\leq t\leq T}}$\footnote{For all $0\leq s\leq t\leq T$, $U(t,s)$ is defined as the evolution operator associated with $\mathcal{A}(t)$, \textit{i.e.}
\begin{equation*}
\frac{\partial}{\partial t}U(t,s) = \mathcal{A}(t) U(t,s),\quad U(s,s)=\mathrm{Id},
\end{equation*}
so that the solution to $\partial_t Z(t) = \mathcal{A}(t) Z(t)$, $Z(s)=Z_s$, is of the form $Z(t) = U(t,s)Z_s$.
}. Indeed, we write $Z(s)=\sum_{k\geq1}Z_k(s)e_k(x)$ with $Z_k(s)\in\mathbb R^n$. Since ${\widetilde A}(t),{\widetilde D}(t)\in M_n(\mathbb R)$ do not depend on $x$, and $\Delta e_k=-\lambda_k e_k$, the operator $\mathcal A(t)$ acts diagonally on each eigenspace $\mathrm{span}(e_k)\otimes\mathbb R^n$, with no coupling between distinct modes:
\begin{equation*}
    U(t,s)Z(s)=\sum_{k=1}^{\infty}U_{\lambda_k}^*(t,s)Z_k(s)\,e_k,
\end{equation*}
where $Z_k(t):=U_{\lambda_k}^*(t,s)Z_k(s)$ solves the finite-dimensional linear ODEs \eqref{eq:decomposition}.
By uniqueness of solutions to linear ODEs, $Z_k(s)=0\ \Rightarrow\ Z_k(t)=0$ for all $t\geq s$. Hence:
\begin{itemize}
\item if $Z(s)\in E_J$ (\textit{i.e.} $Z_k(s)=0$ for $k>J$), then $Z_k(t)=0$ for $k>J$ and all $t\geq s$, so $U(t,s)Z(s)\in E_J$;
\item if $Z(s)\in F_J$ (\textit{i.e.} $Z_k(s)=0$ for $k\leq J$), then $Z_k(t)=0$ for $k\leq J$ and all $t\geq s$, so $U(t,s)Z(s)\in F_J$.
\end{itemize}
Therefore
\begin{equation}\label{stability_E_J}
    U(t,s)E_J\subset E_J,\qquad U(t,s)F_J\subset F_J,\qquad\forall\,{0\leq s\leq t\leq T}.
\end{equation}
We now establish the low-frequency ingredient of the Lebeau--Robbiano
strategy: an observability estimate for each mode system \eqref{eq:decomposition}
on an arbitrary subinterval of ${[\widetilde a,\widetilde b]}$ defined in \eqref{def:atilde}, whose constant is explicit in the
length $\varepsilon$ of the subinterval and grows at most polynomially in
$\lambda_k$. 
\begin{lmm}\label{observabilityineq}
Assume that Assumption~\ref{hyp:USM} holds. Then there exist $C_0>0$ and
$M\in\xN$ such that, for any $\varepsilon\in(0,{\widetilde b-\widetilde a})$, any subinterval
$I=[c,c+\varepsilon]$ of ${[\widetilde a,\widetilde b]}$ of length $\varepsilon$, any $k\ge1$ and
any solution $Z_k$ of the ODE in \eqref{eq:decomposition},
\begin{equation}\label{ine:obs1}
\|Z_k(c+\varepsilon)\|^2\le C(\varepsilon,k)\int_I\|{\widetilde B}^*(t)Z_k(t)\|^2\,dt,
\end{equation}
where
\begin{equation}\label{ine:obsconstant}
C(\varepsilon,k):=C_0\,\varepsilon^{-(2p+1)}(1+\lambda_k)^{2M}.
\end{equation}
\end{lmm}
\begin{proof}
We consider the control problem
dual to \eqref{eq:decomposition}, namely (see \eqref{eq:dualcontrol}, with $\lambda$ replaced by $\lambda_k$)
\begin{equation}\label{eq:dualdecomposition}
\dot y(t)=-{\widetilde A}_{\lambda_k}(t)\,y(t)+{\widetilde B}(t)\,u_k(t),
\end{equation}
where $u_k\in L^2(I;\xR^m)$, and we denote by $C_i^{\lambda_k}$, $i\ge0$, its
Silverman--Meadows matrices \eqref{eq:SMdual}. Let $\varepsilon\in(0,{\widetilde b-\widetilde a})$,
and consider any subinterval $I=[c,c+\varepsilon]$ of ${[\widetilde a,\widetilde b]}$ of length
$\varepsilon$. Now, we follow the construction given in
\cite[page 14]{MR2302744}. From \eqref{eq:USMdual}, we deduce that
\begin{equation}\label{eq:Ginvertible}
\begin{aligned}
   & G_{\lambda_k}(t):=\sum_{i=0}^pC_i^{\lambda_k}(t)\,C_i^{\lambda_k}(t)^*
\;\text{is invertible for every } t\in{[\widetilde a,\widetilde b]},\\&
\;\text{with } \det G_{\lambda_k}(t)=P_p(\lambda_k,{T-t})\ge c_G.
\end{aligned}
\end{equation}
We define, for $j\in\{0,\ldots,p\}$,
$Q_{\lambda_k,j}\in C^\infty([\widetilde a,\widetilde b];M_{m,n}(\xR))$ by
\begin{equation}\label{def:Q_kj}
\begin{aligned}
    &Q_{\lambda_k,j}(t)=C_j^{\lambda_k}(t)^*\,G_{\lambda_k}(t)^{-1},\;\forall t\in [\widetilde a,\widetilde b],
\qquad\text{so that}\\&
\sum_{j=0}^pC_j^{\lambda_k}(t)\,Q_{\lambda_k,j}(t)=I_n .
\end{aligned}
\end{equation}
Let us consider the smaller interval
\begin{gather}\label{eq:Jdefine}
\mathcal J:=\Bigl[c+\frac\varepsilon4,\,c+\frac{3\varepsilon}4\Bigr].
\end{gather}
Fix $\Phi_0\in C^\infty(\xR)$ such that $\Phi_0=0$ on
$(-\infty,\tfrac14]$ and $\Phi_0=1$ on $[\tfrac34,+\infty)$, and set
\begin{equation}\label{def:var-Phi}
\Phi(t):=\Phi_0\Bigl(\frac{t-c}{\varepsilon}\Bigr),\qquad t\in I .
\end{equation}
Clearly, $\Phi\in C^\infty(I)$. By construction,
\begin{gather}\label{eq:Phi=0}
\Phi(c)=0,\;\;\Phi(c+\varepsilon)=1,
\end{gather}
all the derivatives of $\Phi$ vanish on $I\setminus\mathcal J$, and, by scaling,
\begin{equation}\label{lanorme}
\|\Phi^{(i)}\|_\infty\le\frac{\|\Phi_0^{(i)}\|_\infty}{\varepsilon^{i}},
\qquad i\in\{0,\ldots,p+1\}.
\end{equation}
Let $y^1\in\xR^n$, and let $\gamma_k\in C^\infty(I;\xR^n)$ be the solution of
the (backward) Cauchy problem
\begin{equation}\label{eq:gamma1k}
\dot\gamma_k(t)=-{\widetilde A}_{\lambda_k}(t)\gamma_k(t),\;\;\gamma_k(c+\varepsilon)=y^1.
\end{equation}
Let $\Lambda_k\in C^\infty(I;\xR^n)$ be defined by
\begin{equation}\label{eq:Lambdak}
\Lambda_k(t):=\Phi(t)\gamma_k(t),\;\;\forall t\in I .
\end{equation}
From \eqref{eq:Phi=0} and \eqref{eq:gamma1k} one has
\begin{equation}
\Lambda_k(c)=0,\;\;\Lambda_k(c+\varepsilon)=y^1 .
\end{equation}
Let $q_k\in C^\infty(I;\xR^n)$ be defined by
\begin{equation} \label{def:qk}
q_k(t):=\dot\Lambda_k(t)+{\widetilde A}_{\lambda_k}(t)\Lambda_k(t),\;\;\forall t\in I .
\end{equation}
From \eqref{eq:Jdefine}, \eqref{def:var-Phi}, \eqref{eq:gamma1k} and
\eqref{eq:Lambdak}, one has
\begin{equation}\label{qk}
q_k(t)=\dot\Phi(t)\gamma_k(t)\;\;\forall t\in I,
\qquad\text{and}\qquad q_k=0\;\text{on } I\setminus\mathcal J .
\end{equation}
Following \cite[page 14]{MR2302744}, the control $u_k$ steering $0$ to $y^1$
on $I$ is defined by
\begin{equation}\label{def:Uk}
u_k(t)=\sum_{j=0}^p(-1)^j\Bigl(Q_{\lambda_k,j}(t)q_k(t)\Bigr)^{(j)},\;\forall t\in I ,
\end{equation}
so that $u_k=0$ on $I\setminus\mathcal J$ by \eqref{qk}. Let us check that
$u_k$ indeed steers $0$ to $y^1$. Writing $L_k:=\frac{d}{dt}+{\widetilde A}_{\lambda_k}$, the
recursion \eqref{eq:SMdual} reads $C_i^{\lambda_k}=L_kC_{i-1}^{\lambda_k}$, hence
$L_k(C_{i-1}^{\lambda_k}w)=C_i^{\lambda_k}w+C_{i-1}^{\lambda_k}\dot w$ for every
smooth $w$, and by induction
\begin{equation*}
C_j^{\lambda_k}w
=L_k\Bigl(\sum_{i=0}^{j-1}(-1)^iC_{j-1-i}^{\lambda_k}w^{(i)}\Bigr)+(-1)^jC_0^{\lambda_k}w^{(j)} .
\end{equation*}
Applying this with $w=Q_{\lambda_k,j}q_k$, summing over $j$ and using
\eqref{def:Q_kj}, we get $q_k=LR_k+{\widetilde B}u_k$ with
$R_k:=\sum_{j=0}^p\sum_{i=0}^{j-1}(-1)^iC_{j-1-i}^{\lambda_k}(Q_{\lambda_k,j}q_k)^{(i)}$,
and $R_k(c)=R_k(c+\varepsilon)=0$ since $q_k$ vanishes near $c$ and
$c+\varepsilon$ by \eqref{qk}. Therefore, $y:=\Lambda_k-R_k$ solves
\eqref{eq:dualdecomposition} with $u_k$ defined in \eqref{def:Uk}, $y(c)=0$ and $y(c+\varepsilon)=y^1$.
\par We now establish an upper bound for the \(L^2\)-norm of $u_k$. From \eqref{qk}, \eqref{def:Uk} and using the Leibniz rule, we have, $\forall t\in I$,
\begin{equation*}
u_k(t)=\sum_{j=0}^p(-1)^j\sum_{i=0}^j\binom ji\,Q_{\lambda_k,j}^{(j-i)}(t)\,q_k^{(i)}(t),
\qquad
q_k^{(i)}(t)=\sum_{\ell=0}^i\binom i\ell\,\Phi^{(\ell+1)}(t)\,\gamma_k^{(i-\ell)}(t).
\end{equation*}
Using the coercivity assumption \eqref{coercivityofD}, the solution $\gamma_k$
is uniformly bounded with respect to $k$. Indeed,
\begin{equation*}
\frac12\frac{d}{dt}\|\gamma_k(t)\|^2
=-\langle {\widetilde A}_{\lambda_k}(t)\gamma_k(t),\gamma_k(t)\rangle
\ge\alpha\lambda_k\|\gamma_k(t)\|^2-\|{\widetilde A}\|_{L^\infty}\|\gamma_k(t)\|^2
\ge-\|{\widetilde A}\|_{L^\infty}\|\gamma_k(t)\|^2 .
\end{equation*}
Hence, by Gronwall's lemma applied backward from $t=c+\varepsilon$,
\begin{equation*}
\|\gamma_k(t)\|\le C\|y^1\|,\;\;t\in I,
\end{equation*}
with $C$ independent of $k$, $c$ and $\varepsilon$ (but depending on $\tilde a$ and $\tilde b$). Since
$t\mapsto {\widetilde A}_{\lambda_k}(t)$ depends smoothly on $t$ and polynomially on
$\lambda_k$, differentiating \eqref{eq:gamma1k} shows by induction that, for
every $r\in\{0,\ldots,p\}$,
\begin{equation*}
\|\gamma_k^{(r)}(t)\|\le C(1+\lambda_k)^r\|y^1\|,\;\;t\in I .
\end{equation*}
Consequently, by \eqref{lanorme} and \eqref{def:qk},
\begin{equation}\label{eq:qbound}
\|q_k^{(i)}\|_{L^\infty(I)}\le C\,\varepsilon^{-i-1}(1+\lambda_k)^{i}\|y^1\|,
\;\;i=0,\ldots,p .
\end{equation}
Let
\begin{equation*}
h(\lambda_k):=\max_{\substack{0\le j\le p\\0\le r\le p}}\ \sup_{t\in{[\widetilde a,\widetilde b]}}
\Bigl\|\frac{d^r}{dt^r}Q_{\lambda_k,j}(t)\Bigr\| .
\end{equation*}
We now claim that $h(\lambda_k)$ has at most polynomial growth with respect to
$\lambda_k$. From \eqref{eq:SMdual}, the entries of the matrices $C_i^{\lambda_k}(t)$ are polynomials in
$\lambda_k$ of degree at most $i$, with smooth coefficients depending on $t$.
Therefore the entries of $G_{\lambda_k}(t)$ and of $\operatorname{adj}G_{\lambda_k}(t)$
are also polynomials in $\lambda_k$, with smooth coefficients in $t$. Since
$G_{\lambda_k}^{-1}=\operatorname{adj}G_{\lambda_k}/\det G_{\lambda_k}$ and
$\det G_{\lambda_k}\ge c_G$ on ${[\widetilde a,\widetilde b]}$ by \eqref{eq:Ginvertible}, 
and since
$(G_{\lambda_k}^{-1})'
=-G_{\lambda_k}^{-1}G_{\lambda_k}'G_{\lambda_k}^{-1}$,
{a straightforward induction on $r$, based on repeated
differentiation of this identity, shows that, for every
$r\in\{0,\ldots,p\}$, $(G_{\lambda_k}^{-1})^{(r)}$ has at most
polynomial growth in $\lambda_k$, uniformly for
$t\in[\widetilde a,\widetilde b]$.} Hence the entries of $Q_{\lambda_k,j}$ and of their time derivatives up to order $p$
are bounded, uniformly for $t\in[\widetilde a,\widetilde b]$, by polynomials in $\lambda_k$. Thus
there exist $C>0$ and $M_1\in\xN$, independent of $k$, such that
\begin{equation*}
h(\lambda_k)\le C(1+\lambda_k)^{M_1}.
\end{equation*}
Combining this estimate with \eqref{eq:qbound}, and using $\varepsilon\le \min(1,{\widetilde b-\widetilde a})$
to absorb the lower powers of $\varepsilon^{-1}$, we obtain
\begin{equation*}
\|u_k\|_{L^\infty(I)}\le C\,\varepsilon^{-(p+1)}(1+\lambda_k)^{M}\|y^1\|,
\qquad M:=M_1+p .
\end{equation*}
Therefore,
\begin{equation}\label{eq:Ubound}
\|u_k\|_{L^2(I)}\le\sqrt\varepsilon\,\|u_k\|_{L^\infty(I)}
\le C\,\varepsilon^{-p-\frac12}(1+\lambda_k)^{M}\|y^1\| .
\end{equation}
Finally, let $Z_k$ be any solution of the ODE in \eqref{eq:decomposition} on
$I$, and apply the above construction with $y^1:=Z_k(c+\varepsilon)$. Since
$y:=\Lambda_k-R_k$ solves \eqref{eq:dualdecomposition} with $u_k$ defined in \eqref{def:Uk}, $y(c)=0$, $y(c+\varepsilon)=y^1$ and $Z_k'={\widetilde A}_{\lambda_k}^*Z_k$, one has by duality that
$$\frac{d}{dt}\langle y(t),Z_k(t)\rangle=\langle {-\widetilde A}_{\lambda_k}(t)\,y(t)+{\widetilde B}(t)\,u_k(t),Z_k(t)\rangle+\langle y(t), {\widetilde A}_{\lambda_k}^*(t)Z_k(t) \rangle=\langle u_k(t),{\widetilde B}^*Z_k(t)\rangle.$$ Therefore, integrating over $I$,
\begin{equation*}
\|Z_k(c+\varepsilon)\|^2
=\langle y(c+\varepsilon),Z_k(c+\varepsilon)\rangle-\langle y(c),Z_k(c)\rangle
=\int_I\langle u_k(t),{\widetilde B}^*(t)Z_k(t)\rangle\,dt
\le\|u_k\|_{L^2(I)}\,\|{\widetilde B}^*Z_k\|_{L^2(I)} .
\end{equation*}
By \eqref{eq:Ubound} with $\|y^1\|=\|Z_k(c+\varepsilon)\|$, this yields
\begin{equation*}
\|Z_k(c+\varepsilon)\|\le C\,\varepsilon^{-p-\frac12}(1+\lambda_k)^{M}\,\|{\widetilde B}^*Z_k\|_{L^2(I)},
\end{equation*}
which is \eqref{ine:obs1}-\eqref{ine:obsconstant} with $C_0:=C^2$.
\end{proof}

The next lemma is the counterpart of \cite[Lemma~2.1]{MR2679651} in our
non-autonomous setting (see Remark~\ref{rem:strategy}): combining
Lemma~\ref{observabilityineq} with the Lebeau--Robbiano spectral inequality
for the low frequencies and the parabolic dissipation for the high ones, it
provides a local estimate on every subinterval $[c,c+\varepsilon]\subset{[\widetilde a,\widetilde b]}$,
which is then propagated by a telescoping argument.
\begin{lmm}\label{MillerLemma}
There exist $q\in(0,1)$ and a function
$f:(0,{\widetilde b-\widetilde a})\to(0,+\infty)$, with
\begin{equation*}
    f(\varepsilon)\longrightarrow 0
\;\;\text{as }\varepsilon\to0^+,
\end{equation*}
such that, for every $\varepsilon>0$ small enough, every interval $[c,c+\varepsilon]\subset {[\widetilde a,\widetilde b]}$ and every solution
$Z$ of \eqref{eq:obs}, one has
\begin{equation}\label{ine:MillerStep}
f(\varepsilon)\|Z(c+\varepsilon)\|_{L^2(\Omega)^n}^2
-
f(q\varepsilon)\|Z(c)\|_{L^2(\Omega)^n}^2
\le
\int_c^{c+\varepsilon}
\|W^*{\widetilde B}^*(t)Z(t)\|_{L^2(\Omega)^m}^2\,dt.
\end{equation}
Consequently,
\begin{equation}\label{ine:obscost}
\|Z(c+\varepsilon)\|_{L^2(\Omega)^n}^2
\le
\frac{1}{f((1-q)\varepsilon)}
\int_c^{c+\varepsilon}
\|W^*{\widetilde B}^*(t)Z(t)\|_{L^2(\Omega)^m}^2\,dt.
\end{equation}
\end{lmm}

\begin{proof}[Proof of Lemma~\ref{MillerLemma}]
We follow the ideas of \cite{MR2679651}. Let $\varepsilon\in(0,{\widetilde b-\widetilde a})$ and let $[c,c+\varepsilon]\subset{[\widetilde a,\widetilde b]}$. We first prove
\eqref{ine:MillerStep}.

Let $N\in\xN^*$, to be chosen later, and decompose
\begin{equation}\label{def:Z_decomposition}
    Z=Z_{\le N}+Z_{>N},
\end{equation}
where
\begin{equation}
    Z_{\le N}(t,x)=\sum_{k\le N}Z_k(t)e_k(x),
\;\;
Z_{>N}(t,x)=\sum_{k>N}Z_k(t)e_k(x).
\end{equation}
By \eqref{stability_E_J}, $Z_{\le N}$ and $Z_{>N}$ are the solutions of \eqref{eq:obs} with  data at time $t=c$ given by the orthogonal projections of $Z(c)$ on $E_N$ and $F_N$. We apply Lemma~\ref{observabilityineq} on the interval
\begin{equation*}
    I_1:=\left[c+\frac{\varepsilon}{2},c+\varepsilon\right].
\end{equation*}
Taking
$c'=c+\frac{\varepsilon}{2}$ and $\delta=\frac{\varepsilon}{2}$ in Lemma~\ref{observabilityineq}
\eqref{ine:obs1} yields
\begin{equation}\label{ine:obs2proof}
 \Vert Z_k(c+\varepsilon) \Vert^2
 \leq C\left(\frac{\varepsilon}{2},k\right)
 \int_{c+\frac{\varepsilon}{2}}^{c+\varepsilon}
 \Vert {\widetilde B}^*(t)Z_k(t)\Vert^2\,dt.
\end{equation}
Moreover, by the expression \eqref{ine:obsconstant} of the constant,
\begin{equation}\label{ine:obsconstant-proof}
C\Bigl(\frac\varepsilon2,k\Bigr)=2^{2p+1}\,C(\varepsilon,k)
=2^{2p+1}C_0\,\varepsilon^{-(2p+1)}(1+\lambda_k)^{2M},
\end{equation}
with $C_0>0$ a constant, independent of $k$, $N$,
$c$ and $\varepsilon$. Summing \eqref{ine:obs2proof} on $I_1$ over $k\le N$ and using
\eqref{ine:obsconstant-proof}, we obtain
\begin{equation} \label{eq:lowfreq}
\sum_{k\le N}\|Z_k(c+\varepsilon)\|^2
\le C\,\varepsilon^{-(2p+1)}(1+\lambda_N)^{2M}
\int_{I_1}\sum_{k\le N}\|{\widetilde B}^*(t)Z_k(t)\|^2\,dt .
\end{equation}
To pass from the full observation ${\widetilde B}^*Z_k$ to the localized one, we use the
Lebeau--Robbiano spectral inequality, in the form
given in \cite{JerisonLebeau1999,LebeauZuazua1998}: for every nonempty open
$\omega_i\subset\Omega$ there exist $C_i,c_i>0$ such that, for every $N\ge1$
and every $(a_k)_{k\le N}\in\xR^N$,
\begin{equation}\label{ine:LR}
\sum_{k\le N}|a_k|^2\le C_ie^{c_i\sqrt{\lambda_N}}\int_{\omega_i}\Bigl|\sum_{k\le N}a_ke_k(x)\Bigr|^2dx .
\end{equation}
Applying \eqref{ine:LR} for each fixed $t\in I_1$ with $a_k=({\widetilde B}^*(t)Z_k(t))_j$,
$1\le j\le m$, summing over $j$ and integrating over $I_1$ yields
\begin{equation}\label{ine:LRtime}
\int_{I_1}\|{\widetilde B}^*(t)Z_{\le N}(t)\|_{L^2(\Omega)^m}^2dt
\le Ce^{C\sqrt{\lambda_N}}\int_{I_1}\|W^*{\widetilde B}^*(t)Z_{\le N}(t)\|_{L^2(\Omega)^m}^2dt,
\end{equation}
where we used $\|{\widetilde B}^*(t)Z_{\le N}(t)\|_{L^2(\Omega)^m}^2=\sum_{k\le N}\|{\widetilde B}^*(t)Z_k(t)\|^2$
by orthonormality of $(e_k)$.
From \eqref{eq:lowfreq} and using  $Z_{\le N}=Z-Z_{>N}$, this gives
\begin{align}
\|Z_{\le N}(c+\varepsilon)\|_{L^2(\Omega)^n}^2
&\le
\frac{C}{\varepsilon^{2p+1}}
(1+\lambda_N)^{2M} e^{C\sqrt{\lambda_N}}
\int_{c}^{c+\varepsilon}
\|W^*{\widetilde B}^*(t)Z(t)\|_{L^2(\Omega)^m}^2\,dt
\nonumber\\
&\;+
\frac{C}{\varepsilon^{2p+1}}
(1+\lambda_N)^{2M} e^{C\sqrt{\lambda_N}}
\int_{I_1}
\|W^*{\widetilde B}^*(t)Z_{>N}(t)\|_{L^2(\Omega)^m}^2\,dt.
\label{ineq:low}
\end{align}

We now estimate the high-frequency part. By the coercivity assumption on ${\widetilde D}$ (which is an immediate consequence of \eqref{coercivityofD}), there
exist constants $C>0$ and $\mu>0$ such that, for every $k$ large enough,
\begin{equation*}
    \|Z_k(t)\|^2
\le
C e^{-\mu\lambda_k(t-c)}\|Z_k(c)\|^2,
\;\; t\in[c,c+\varepsilon].
\end{equation*}
Since $t-c\leq \varepsilon < {\widetilde b-\widetilde a}$, increasing $C$ if necessary, this estimate may be assumed to hold for every $k>N$.
Therefore,
\begin{equation}\label{ine:dissipation}
    \|Z_{>N}(c+\varepsilon)\|_{L^2(\Omega)^n}^2
\le
C e^{-\mu\lambda_{N+1}\varepsilon}
\|Z(c)\|_{L^2(\Omega)^n}^2,
\end{equation}
and, since $t-c\ge\varepsilon/2$ on $I_1$, and by the analogue of \eqref{eq:admissibility} with initial time $c+\varepsilon/2$,

\begin{align}\label{ineq:int}
    \int_{I_1}
\|W^*{\widetilde B}^*(t)Z_{>N}(t)\|_{L^2(\Omega)^m}^2\,dt &\leq C\|Z_{>N}(c+\varepsilon/2)\|_{L^2(\Omega)^n}^2\nonumber\\
&\le 
C e^{-\mu\lambda_{N+1}\varepsilon/2}
\|Z(c)\|_{L^2(\Omega)^n}^2.
\end{align}
Recall that, thanks to \eqref{def:direct_sum} and  \eqref{def:Z_decomposition}, we have
\begin{equation*}
    \|Z(c+\varepsilon)\|_{L^2(\Omega)^n}^2
    =
    \|Z_{\leq N}(c+\varepsilon)\|_{L^2(\Omega)^n}^2
    +
    \|Z_{>N}(c+\varepsilon)\|_{L^2(\Omega)^n}^2.
\end{equation*}
Combining this identity with the estimates \eqref{ineq:low}, 
\eqref{ine:dissipation} and \eqref{ineq:int}, we obtain
\begin{align}
\|Z(c+\varepsilon)\|_{L^2(\Omega)^n}^2
&\le
\frac{C}{\varepsilon^{2p+1}}
(1+\lambda_N)^{2M} e^{C\sqrt{\lambda_N}}
\int_{c}^{c+\varepsilon}
\|W^*{\widetilde B}^*(t)Z(t)\|_{L^2(\Omega)^m}^2\,dt
\nonumber\\
&\;+
C\left[
1+
\frac{1}{\varepsilon^{2p+1}}
(1+\lambda_N)^{2M} e^{C\sqrt{\lambda_N}}
\right]
e^{-\mu\lambda_{N+1}\varepsilon/2}
\|Z(c)\|_{L^2(\Omega)^n}^2.
\label{ineq:one-step-pre}
\end{align}

{
We now choose the spectral threshold. Let $L_0>0$, to be fixed later, and set
\begin{equation}\label{def:threshold}
    \Lambda(\varepsilon):=\lambda_1+\frac{L_0}{\varepsilon^2}.
\end{equation}
Choose $N=N(\varepsilon)$ such that
\begin{equation}\label{ineq:threshold}
    \lambda_N\le \Lambda(\varepsilon)<\lambda_{N+1}.
\end{equation}
Fix $C\ge1$ so that \eqref{ineq:one-step-pre} holds. For
$0<\varepsilon\le1$, we have
\[
\sqrt{\lambda_N}\le\sqrt{\lambda_1}+\frac{\sqrt{L_0}}{\varepsilon},
\qquad
1+\lambda_N\le(1+\lambda_1+L_0)\varepsilon^{-2}.
\]
Thus, with $\theta_0:=2p+1+4M$, there exists $C_{L_0}\ge1$,
independent of $c$ and $\varepsilon$, such that
\begin{equation*}
C\left[1+\frac{1}{\varepsilon^{2p+1}}(1+\lambda_N)^{2M}
e^{C\sqrt{\lambda_N}}\right]
\le C_{L_0}\varepsilon^{-\theta_0}
e^{C\sqrt{L_0}/\varepsilon}.
\end{equation*}
For fixed $L_0$, since
$\varepsilon\log(C_{L_0}\varepsilon^{-\theta_0})\to0$ as
$\varepsilon\to0^+$, there exists
$\varepsilon_0\in(0,\min\{1,\widetilde b-\widetilde a\})$ such that
$C_{L_0}\varepsilon^{-\theta_0}\le e^{1/\varepsilon}$ for every
$0<\varepsilon\le\varepsilon_0$. Setting $K_0:=C\sqrt{L_0}+1$, we obtain
\begin{equation}\label{ineq:term1-final}
\frac{C}{\varepsilon^{2p+1}}(1+\lambda_N)^{2M}
e^{C\sqrt{\lambda_N}}\le e^{K_0/\varepsilon},
\end{equation}
and
\begin{equation}\label{ineq:bracket-final}
C\left[1+\frac{1}{\varepsilon^{2p+1}}(1+\lambda_N)^{2M}
e^{C\sqrt{\lambda_N}}\right]\le e^{K_0/\varepsilon},
\end{equation}
for every $0<\varepsilon\le\varepsilon_0$.
By \eqref{ineq:threshold},
\begin{equation}\label{ine:decay_exp}
e^{-\mu\lambda_{N+1}\varepsilon/2}
\le e^{-\mu L_0/(2\varepsilon)}.
\end{equation}
Choose $L_0$ large enough that
$C_1:=\mu L_0/2-K_0>0$, and fix the corresponding $\varepsilon_0$.
This is possible since $K_0=C\sqrt{L_0}+1$. Combining
\eqref{ineq:bracket-final} and \eqref{ine:decay_exp} gives
\begin{equation}\label{ineq:term2-small}
C\left[1+\frac{1}{\varepsilon^{2p+1}}(1+\lambda_N)^{2M}
e^{C\sqrt{\lambda_N}}\right]
e^{-\mu\lambda_{N+1}\varepsilon/2}
\le e^{-C_1/\varepsilon}.
\end{equation}
Together with \eqref{ineq:one-step-pre} and \eqref{ineq:term1-final}, this yields
\begin{equation} \label{ine:all}
\begin{split}
\|Z(c+\varepsilon)\|_{L^2(\Omega)^n}^2
\le{}& e^{K_0/\varepsilon}
\int_c^{c+\varepsilon}
\|W^*\widetilde B^*(t)Z(t)\|_{L^2(\Omega)^m}^2\,dt\\
&+e^{-C_1/\varepsilon}\|Z(c)\|_{L^2(\Omega)^n}^2,
\end{split}
\end{equation}
for every $0<\varepsilon\le\varepsilon_0$.
Set $f(\varepsilon):=e^{-K_0/\varepsilon}$ and choose $q\in(0,1)$ such that
$K_0(q^{-1}-1)\le C_1$. Then
$f(\varepsilon)e^{K_0/\varepsilon}=1$ and
$f(\varepsilon)e^{-C_1/\varepsilon}\le f(q\varepsilon)$.
Multiplying \eqref{ine:all} by $f(\varepsilon)$ gives
\eqref{ine:MillerStep} for $\varepsilon\le\varepsilon_0$.

}

It remains to telescope. Fix $[c,c+\varepsilon]\subset{[\widetilde a,\widetilde b]}$ with $\varepsilon\le\varepsilon_0$ and define
\begin{equation*}
    t_\ell:=c+q^\ell\varepsilon,
\;\;
\delta_\ell:=t_\ell-t_{\ell+1}
=
(1-q)q^\ell\varepsilon,
\;\; \ell\in\xN.
\end{equation*}
Then $t_0=c+\varepsilon$, $t_\ell\to c$ as $l \to +\infty$, and $q\delta_\ell=\delta_{\ell+1}.$
Applying \eqref{ine:MillerStep} on each interval $[t_{\ell+1},t_\ell]$, whose length is
$\delta_\ell$, gives
\begin{equation*}
    f(\delta_\ell)\|Z(t_\ell)\|_{L^2(\Omega)^n}^2
-
f(\delta_{\ell+1})\|Z(t_{\ell+1})\|_{L^2(\Omega)^n}^2
\le
\int_{t_{\ell+1}}^{t_\ell}
\|W^*{\widetilde B}^*(t)Z(t)\|_{L^2(\Omega)^m}^2\,dt.
\end{equation*}
Summing from $\ell=0$ to $L$, we obtain
\begin{equation*}
    f(\delta_0)\|Z(c+\varepsilon)\|_{L^2(\Omega)^n}^2 -
f(\delta_{L+1})\|Z(t_{L+1})\|_{L^2(\Omega)^n}^2
\le
\int_{t_{L+1}}^{c+\varepsilon}
\|W^*{\widetilde B}^*(t)Z(t)\|_{L^2(\Omega)^m}^2\,dt.
\end{equation*}
Letting $L\to+\infty$, using $t_{L+1}\to c$ and $f(\delta_{L+1})\to0$, we get
\begin{equation*}
    f((1-q)\varepsilon)\|Z(c+\varepsilon)\|_{L^2(\Omega)^n}^2
\le
\int_c^{c+\varepsilon}
\|W^*{\widetilde B}^*(t)Z(t)\|_{L^2(\Omega)^m}^2\,dt,
\end{equation*}
which is \eqref{ine:obscost} for $\varepsilon\le\varepsilon_0$.
\end{proof}
\begin{proof}[Proof of Theorem~\ref{main:thm}]
Let $[c,c+\varepsilon]\subset{[\widetilde a,\widetilde b]}$ with $\varepsilon\le\varepsilon_0$, where
$\varepsilon_0$ is given by Lemma~\ref{MillerLemma}. By \eqref{ine:obscost},
\[
\|Z(c+\varepsilon)\|_{L^2(\Omega)^n}^2\le\frac{1}{f((1-q)\varepsilon)}
\int_c^{c+\varepsilon}\|W^*{\widetilde B}^*(t)Z(t)\|_{L^2(\Omega)^m}^2\,dt .
\]
Moreover, by the energy estimate for \eqref{eq:obs},
$\|Z(T)\|_{L^2(\Omega)^n}^2\le C\|Z(c+\varepsilon)\|_{L^2(\Omega)^n}^2$. Since $[c,c+\varepsilon]\subset [0,T]$, we deduce that
\eqref{ine:obs} holds, which concludes the proof.
\end{proof}

\subsection{Sufficient conditions for Assumption~\ref{hyp:USM}}
\label{sec:sufficient}

All the sufficient conditions of this section rely on the following
pointwise condition, which is the direct spectral analogue of the
Silverman--Meadows rank condition at a fixed time.

\begin{hyp}[Pointwise spectral Silverman--Meadows condition]\label{hyp:pointwise}
There exists $\bar t\in[0,T]$ such that
\begin{equation}\label{APO}
\operatorname{Span}\Big\{B_i^{\lambda_k}(\bar t)U;\ U\in\xR^m,\ i\in\xN\Big\}=\xR^n,
\qquad\forall k\in\xN^*.
\end{equation}
\end{hyp}

The equality \eqref{APO} means that for every $k$ there exists
$p_k\in\xN$ such that $\operatorname{rank}\mathcal K_{\lambda_k,p_k}(\bar t)=n$,
\textit{i.e.} $P_{p_k}(\lambda_k,\bar t)>0$, where $p_k$ may \textit{a priori} depend
on $k$. The following lemma proves that, in fact, $p_k$ can be chosen independently of $k$. 

\begingroup

\begin{lmm}\label{indep}
Assume that the pointwise spectral Silverman--Meadows condition \eqref{APO}
holds. Then there exists $p\in\xN$,
independent of $k$, such that
\begin{equation}\label{hyp:P_p>0}
P_p(\lambda_k,\bar t)>0,\qquad\forall k\in\xN^*.
\end{equation}
\end{lmm}

\begin{proof}
By the pointwise spectral Silverman--Meadows condition \eqref{APO},
there exists $p_1\in\xN$ such that
$P_{p_1}(\lambda_1,\bar t)>0$. Hence
\[
Q(\lambda):=P_{p_1}(\lambda,\bar t)
\]
is a nonzero polynomial. Since it has finitely many real
roots and $\lambda_k\to+\infty$, the set
\[
F:=\{k\in\xN^*:Q(\lambda_k)=0\}
\]
is finite. For each $k\in F$, the pointwise spectral Silverman--Meadows condition \eqref{APO}
provides an order $p_k\in\xN$ such that
$P_{p_k}(\lambda_k,\bar t)>0$. Set
\[
p:=\max\bigl(\{p_1\}\cup\{p_k:k\in F\}\bigr).
\]
For $k\notin F$, we have $P_{p_1}(\lambda_k,\bar t)>0$. Thus, from \eqref{eq:equi_rank},  the matrix $K_{\lambda_k,p_1}(\bar t)$
has rank $n$. For $k\in F$, from \eqref{APO}, the matrix
$K_{\lambda_k,p_k}(\bar t)$ has rank $n$.
Since the rank of $K_{\lambda,s}(\bar t)$ is non-decreasing
in $s$, we obtain
$\operatorname{rank}K_{\lambda_k,p}(\bar t)=n$
for every $k\in\xN^*$, which is equivalent to
\eqref{hyp:P_p>0}.
\end{proof}
\endgroup

For each fixed $k$, under Assumption~\ref{hyp:pointwise}, the Silverman--Meadows criterion
\cite[Theorem~1.18]{MR2302744} then gives the controllability of the $k$-th
mode system on any interval containing $\bar t$. However, \eqref{APO} does
not by itself provide an interval on which \eqref{hyp:USM:rank} holds
uniformly in $k$, as the following example shows.

\begin{rmrk}\label{rem:notpointwise}
Take
$n=3$, $m=1$, $\bar t=T/2$,
\begin{equation*}
B=e_1,\qquad
A=\begin{pmatrix}0&0&0\\-1&0&0\\0&-1&0\end{pmatrix},\qquad
D(t)=I_3-\begin{pmatrix}0&0&0\\ s(t)&0&0\\0&\sigma(t)&0\end{pmatrix},
\end{equation*}
where $s,\sigma\in C^\infty([0,T])$ take values in $[0,\tfrac14]$, so that
$D$ satisfies \eqref{coercivityofD} with $\alpha=\tfrac34$. Consider, for the moment, any $\lambda\in\mathbb R$. A direct computation
gives
\begin{align*}A_\lambda(t)&=
\begin{pmatrix}
-\lambda & 0 & 0\\
-1+\lambda s(t) & -\lambda & 0\\
0 & -1+\lambda\sigma(t) & -\lambda
\end{pmatrix},\qquad
B_1^\lambda=\begin{pmatrix}\lambda\\ 1-\lambda s\\ 0\end{pmatrix},\\
B_2^\lambda &=\begin{pmatrix}\lambda^2\\ -\lambda\dot s+2\lambda(1-\lambda s)\\ (1-\lambda s)(1-\lambda\sigma)\end{pmatrix},
\end{align*}
so that $[B_0^\lambda,B_1^\lambda,B_2^\lambda]$ is upper triangular and
\begin{equation}\label{eq:det-global}
\det\bigl[B_0^\lambda(t),B_1^\lambda(t),B_2^\lambda(t)\bigr]=(1-\lambda s(t))^2\,(1-\lambda\sigma(t)).
\end{equation}
Moreover, if $s\equiv1/\lambda$ on an interval $I$, the entry $(2,1)$ of
$A_\lambda(t)$ vanishes on $I$, the first column of $A_\lambda$ is
$-\lambda e_1$, and by induction $B_i^\lambda(t)=\lambda^ie_1$ for all $i\ge0$
and $t\in I$; if $\sigma\equiv1/\lambda$ on an interval $I$ (not necessarily the same as before), the entry
$(3,2)$ of $A_\lambda(t)$ vanishes, the subspace $\{v_3=0\}$ is invariant by
$A_\lambda(t)$, and all $B_i^\lambda(t)$, $t\in I$, have a
vanishing third component. In both cases
\begin{equation}\label{eq:dead}
\operatorname{rank}\mathcal K_{\lambda,p}(t)\le2<n
\qquad\forall p\in\xN,\ \forall t\in I .
\end{equation}

Let us go back to the original question, which reduces to choosing $s$ and $\sigma$ such that there exists no interval $[a,b]\subset I$ on which \eqref{hyp:USM:rank} holds
uniformly in $k$. {Since the Dirichlet eigenvalues $\lambda_k$ of
$-\Delta$ on $\Omega$ tend to $+\infty$, we can choose an increasing
sequence of indices $(k_j)_{j\ge1}$ such that
$\lambda_{k_j}\ge2^{2^j}$.} Set $\ell_j:=1/\lambda_{k_j}$, $a_j:=T/(2j)$ and
$I_j:=[\bar t+a_{j+1},\bar t+a_j]$, so that $\bigcup_{j\ge1}I_j=(\bar t,T]$.
Fix $\theta\in C^\infty(\mathbb R)$, non-decreasing, with $\theta=0$ on $(-\infty,0]$
and $\theta=1$ on $[1,+\infty)$, and set, for $j\ge2$ and $t\in I_j$,
\begin{equation*}
\tau_j(t):=\ell_{j+1}+(\ell_{j-1}-\ell_{j+1})\,
\theta\Bigl(\frac{t-\bar t-a_{j+1}}{a_j-a_{j+1}}\Bigr),
\end{equation*}
which increases from $\ell_{j+1}$ at the left endpoint of $I_j$ to
$\ell_{j-1}$ at its right endpoint. On $(\bar t,T]$, define
\begin{equation*}
s:=\begin{cases}\ell_j & \text{on } I_j,\ j \text{ odd},\\ \tau_j & \text{on } I_j,\ j \text{ even},\end{cases}
\qquad
\sigma:=\begin{cases}\ell_j & \text{on } I_j,\ j \text{ even},\\ \tau_j & \text{on } I_j,\ j\ge3 \text{ odd},\\ \ell_2 & \text{on } I_1.\end{cases}
\end{equation*} The values match at the endpoints of the
$I_j$ and $\theta$ is flat at $0$ and $1$, so $s,\sigma\in C^\infty((\bar t,T])$. One has, for any  $t\in I_j$, for every $r\ge1$,
\begin{equation*}
|s^{(r)}(t)|\le\|\theta^{(r)}\|_\infty\,\ell_{j-1}\,(a_j-a_{j+1})^{-r}
\le\|\theta^{(r)}\|_\infty\,2^{-2^{j-1}}\Bigl(\frac{2j(j+1)}{T}\Bigr)^{r}
\xrightarrow[j\to\infty]{}0,
\end{equation*}
and $0\le s\le\ell_{j-1}\rightarrow 0$ as $j\rightarrow +\infty$; the same holds for $\sigma$. Hence $s$ and
$\sigma$ extend to $C^\infty$ functions on $[\bar t,T]$, flat at $\bar t$, with
$s(\bar t)=\sigma(\bar t)=0$, and we extend them evenly to $[0,\bar t]$.

With this choice: at $t=\bar t$, $s=\sigma=0$ and, by \eqref{eq:det-global},
$\operatorname{rank}\mathcal K_{\lambda,2}(\bar t)=3$ for every $\lambda\in\xR$,
so that \eqref{APO} holds with $p_k=2$ for all $k$; and every $t\neq\bar t$
belongs to some $I_j$ (or its mirror image) on which $s\equiv1/\lambda_{k_j}$
or $\sigma\equiv1/\lambda_{k_j}$, so that, by \eqref{eq:dead}, the mode $k_j$
violates the rank condition at every order $p$ on its interior. Consequently,
every nontrivial interval $[a,b]\subset[0,T]$ contains a subinterval on which
\eqref{hyp:USM:rank} fails for some $k$, whatever $p$: Assumption~\ref{hyp:USM}
is violated for all $p$ and all $[a,b]$, although \eqref{APO} holds at
$\bar t$.

This example shows that, for smooth coefficients, \eqref{APO} does not imply Assumption~\ref{hyp:USM} without an additional assumption such as
those of Corollaries~\ref{prop:pointwise-stable-degree},
~\ref{cor:analytic-weak} and~\ref{prop:dominant-kalman}. 
Whether the system of this example is nevertheless null-controllable is an open question, which we do not address.
\end{rmrk}

For fixed $p\in\xN$ and $t\in[0,T]$, the function $\lambda\mapsto P_p(\lambda,t)$
is a polynomial (see Remark~\ref{rem:polynomial}) and we denote by
\begin{equation}\label{convention}
d_p(t):=\deg_\lambda P_p(\lambda,t),\qquad\text{with } d_p(t)=-\infty \text{ if } P_p(\cdot,t)\equiv0.
\end{equation}

{The following corollary gives a sufficient condition,
in terms of the local behaviour of $d_p$, for the pointwise rank
condition to imply Assumption~\ref{hyp:USM}.}

\begin{crllr}\label{prop:pointwise-stable-degree}
Assume that Assumption~\ref{hyp:pointwise} holds, and choose $p$ as in \eqref{hyp:P_p>0}. Assume that $d_p$ is locally maximal at $\bar t$. Then Assumption~\ref{hyp:USM}
holds on an interval containing $\bar t$. Consequently, System~\eqref{eq:intro_heat}
is null-controllable at time $T$.
\end{crllr}

\begin{proof}
Let $d:=d_p(\bar t)$. By \eqref{hyp:P_p>0}, the polynomial $P_p(\cdot,\bar t)$
is nonzero, so $d\ge0$ and $\pi_d(\bar t)\neq0$; dividing
$P_p(\lambda_k,\bar t)\ge0$ by $\lambda_k^d$ and letting $k\to\infty$ gives
$\pi_d(\bar t)\ge0$, hence $\pi_d(\bar t)>0$. By local maximality of $d_p$ and
continuity of $\pi_d$, there is a compact neighbourhood $V$ of $\bar t$ in
$[0,T]$ on which $\pi_s\equiv0$ for $s>d$ and $\pi_d>0$. Arguing as in the proof
of Lemma~\ref{lem:uniformGram}, the term $\pi_d(t)\lambda^d$ dominates the others
uniformly on $V$, so there is $\lambda_0$ such that $P_p(\lambda_k,t)>0$ for all
$t\in V$ and all $k$ with $\lambda_k\ge\lambda_0$. For the finitely many $k$
with $\lambda_k<\lambda_0$, \eqref{hyp:P_p>0} and continuity give
$P_p(\lambda_k,\cdot)>0$ on a smaller neighbourhood of $\bar t$. Hence
\eqref{hyp:USM:rank} holds on a nontrivial closed interval containing $\bar t$.
\end{proof}
When the coefficients are analytic in time, \eqref{APO} alone suffices, as we prove in what follows.
\begin{crllr}\label{cor:analytic-weak}
Assume that $A$, $B$ and $D$ are analytic with respect to $t$ and that
Assumption~\ref{hyp:pointwise} holds. Then Assumption~\ref{hyp:USM} holds.
Consequently, System~\eqref{eq:intro_heat} is null-controllable at time $T$.
\end{crllr}

\begin{proof}
We use the following well-known result: a real analytic function on $[0,T]$
which is not identically zero has finitely many zeros. Let $p$ be given by
Lemma~\ref{indep}, so that $P_p(\lambda_k,\bar t)>0$ for every
$k\in\mathbb N^*$. Since $A$, $B$, $D$ are analytic, the recursion
\eqref{def:kalman} shows by induction that the coefficients of each $B_i^\lambda$,
as a polynomial in $\lambda$, are analytic in $t$, and this property is
preserved by the algebraic operations defining $P_p$; hence the coefficients
$\pi_s$ of $P_p(\lambda,t)=\sum_{s=0}^{2np}\pi_s(t)\lambda^s$ are analytic on
$[0,T]$. As $P_p(\lambda_1,\bar t)>0$, they do not all vanish identically, and
we set $d:=\max\{s:\ \pi_s\not\equiv 0\}$. Since $P_p(\lambda,t)\ge 0$ for all
$\lambda\in\mathbb R$, dividing by $\lambda^d$ and letting $\lambda\to+\infty$
gives $\pi_d\ge 0$ on $[0,T]$; being analytic and not identically zero,
$\pi_d$ has finitely many zeros, so we may choose a nontrivial closed interval
$[a,b]\subset[0,T]$ on which $\pi_d>0$. Arguing as in the proof of
Lemma~\ref{lem:uniformGram}, the term $\pi_d(t)\lambda^d$ dominates the others
uniformly on $[a,b]$, so there is $\lambda_0$ such that $P_p(\lambda_k,t)>0$
for all $t\in[a,b]$ and all $k$ with $\lambda_k\ge\lambda_0$. For the finitely
many $k$ with $\lambda_k<\lambda_0$, the function $P_p(\lambda_k,\cdot)$ is
analytic and positive at $\bar t$, hence has finitely many zeros, and we
shrink $[a,b]$ to a nontrivial closed interval avoiding them. Then
$P_p(\lambda_k,t)>0$ for all $t\in[a,b]$ and all $k\in\mathbb N^*$, which is
\eqref{hyp:USM:rank}.
\end{proof}

Finally, since the leading term of $B_i^\lambda$ in $\lambda$ is $\lambda^iD^iB$
(Remark~\ref{rem:polynomial}), a Kalman condition on the pair $(D(\bar t),B(\bar t))$
controls all sufficiently high frequencies near $\bar t$, and \eqref{APO} is then
only needed for the finitely many remaining modes.
\begin{crllr}\label{prop:dominant-kalman}
Assume that Assumption~\ref{hyp:pointwise} holds and that
\begin{equation}\label{dominant_kalman}
\operatorname{rank}\big[B(\bar t),D(\bar t)B(\bar t),\ldots,D(\bar t)^{n-1}B(\bar t)\big]=n .
\end{equation}
Then Assumption~\ref{hyp:USM} holds on an interval containing $\bar t$.
Consequently, System~\eqref{eq:intro_heat} is null-controllable at time $T$.
\end{crllr}

Condition \eqref{dominant_kalman} is a Kalman condition on the dominant part
$\lambda^i D^i B$ of the matrices $B_i^\lambda$ (Remark~\ref{rem:polynomial});
by the Cayley--Hamilton theorem, replacing $n-1$ by any $p\ge n-1$ in
\eqref{dominant_kalman} gives an equivalent condition. It ensures
\eqref{hyp:USM:rank} near $\bar t$ for all large frequencies, while
\eqref{APO} takes care of the finitely many remaining modes.

\begin{proof}
Set $Q(t):=[B(t),D(t)B(t),\ldots,D(t)^{n-1}B(t)]$. By \eqref{dominant_kalman},
there is a set $S$ of $n$ column indices with $\det Q^S(\bar t)\neq0$. Let $d$
be the sum of the block indices $i\in\{0,\ldots,n-1\}$ of the columns in $S$.
By Remark~\ref{rem:polynomial}, the corresponding minor
$\mathrm{MI}_\lambda(t):=\det\mathcal K^S_{\lambda,n-1}(t)$ is a polynomial in
$\lambda$ of degree at most $d$, with leading coefficient $\det Q^S(t)$ and
smooth lower-order coefficients. By continuity, $|\det Q^S|\ge c_0>0$ on a
compact neighbourhood $V$ of $\bar t$ in $[0,T]$, and, arguing as in the
proof of Lemma~\ref{lem:uniformGram}, there is $\lambda_0$ such that
$\mathrm{MI}_{\lambda_k}(t)\neq0$, hence $\operatorname{rank}\mathcal K_{\lambda_k,n-1}(t)=n$,
for all $t\in V$ and all $k$ with $\lambda_k\ge\lambda_0$. For the finitely
many $k$ with $\lambda_k<\lambda_0$, \eqref{APO} provides $p_k$ with
$P_{p_k}(\lambda_k,\bar t)>0$. Set $p:=\max(n-1,\max_{\lambda_k<\lambda_0}p_k)$.
Since the rank of $\mathcal K_{\lambda,p}$ is non-decreasing in $p$,
\eqref{hyp:USM:rank} holds with this $p$ on $V$ for $\lambda_k\ge\lambda_0$,
and, by continuity of $P_p(\lambda_k,\cdot)$, on a smaller neighbourhood of
$\bar t$ for the finitely many other $k$. Hence \eqref{hyp:USM:rank} holds on
a nontrivial closed interval containing $\bar t$.
\end{proof}

\section*{Acknowledgements}
Pierre Lissy is supported by the French Agence Nationale de la Recherche
(Grant ANR-22-CPJ2-0138-01 and Grant ANR-26-CE40-0732-01). Thibault Liard  is  supported  by the French Agence Nationale de la Recherche (Grant ANR-24-CE46-7295).
\vspace{0.3cm}

\footnotesize

\bibliographystyle{plain}
\bibliography{bibliography}

\end{document}